\documentstyle[12pt,amscd]{amsart}
\theoremstyle{plain}
 \newtheorem{thm}{Theorem}[section]
 \newtheorem{lem}[thm]{Lemma}
 \newtheorem{prop}[thm]{Proposition}
 \newtheorem{cor}[thm]{Corollary}
\theoremstyle{definition}
 \newtheorem{defn}{Definition}[section]
\theoremstyle{remark}
 \newtheorem{rem}{Remark}[section]

\newcommand{\Ext}{\operatorname{Ext}}
\newcommand{\Hom}{\operatorname{Hom}}
\newcommand{\codim}{\operatorname{codim}}
\newcommand{\im}{\operatorname{im}}
\newcommand{\Aut}{\operatorname{Aut}}
\newcommand{\rk}{\operatorname{rk}}
\newcommand{\gr}{\operatorname{gr}}

\newcommand{\union}{\operatornamewithlimits{\cup}}

\numberwithin{equation}{section}
\begin{document}

% End-of-Header

\title{Chamber structure of polarizations and\\
the moduli of stable sheaves on a ruled surface}
\author{K\={o}ta Yoshioka\\
 Dept. of Math. Kyoto Univ.\\
 Fax number 75-753-3711\\
 running title\\
 Moduli of stable sheaves on a ruled surface}
 \address{Department of Mathematics, Faculity of Science, Kyoto University,
606-01, Japan }
\thanks{Supported by JSPS Research Fellowships for Young Scientists}
 \maketitle

\section{Introduction}
Let $X$ be a smooth projective surface defined over $\Bbb C$ and $H$ an ample
divisor on $X$.
Let $M_H(r;c_1,c_2)$ be the moduli space of stable sheaves of rank $r$ whose
Chern classes
$(c_1,c_2) \in H^2(X,\Bbb Q)\times H^4(X,\Bbb Q)$ and
$\overline{M}_H(r;c_1,c_2)$ the
Gieseker-Maruyama compactification of $M_H(r;c_1,c_2)$.
When $r=2$, these spaces are extensively studied by many authors.
When $r \geq 3$, Drezet and Le-Potier [D1],[D-L] investigated the structure of
moduli spaces on $\Bbb P^2$, and Rudakov [R] treated moduli spaces on $\Bbb P^1
\times \Bbb P^1$.
In this paper, we shall consider moduli spaces of rank $r \geq 3$ on a ruled
surface which is not rational.
In particular, we shall compute the Picard group of
$\overline{M}_H(r;c_1,c_2)$.
Let $\pi:X \to C$ be the fibration, $f$ a fibre of $\pi$ and $C_0$ a minimal
section of $\pi$ with $(C_0^2)=-e$.
We assume that $e>2g-2$, where $g$ is the genus of $C$.
Then $K_X$ is effective, and hence $(K_X,H)<0$ for any ample divisor $H$.
In particular, $M_H(r;c_1,c_2)$ is smooth with the expected dimension $2r^2
\Delta-r^2(1-g)+1$.

In section 2, we shall generalize the chamber structure of Qin [Q2].
As an application, we shall consider the difference of Betti numbers of moduli
spaces on a ruled surface.
Although we cannot generalize the method in [Y2, 0] directly,
by using Qin's method we can generalize it to any rank case.
In [Y2], we computed the number of $\mu$-semi-stable sheaves of rank 2 on a
ruled surface defind over $\Bbb F_q$.
So, in principle, we can compute the Betti numbers of $M_H(3;c_1,c_2)$ on $\Bbb
P^2$.
Combining chamber structure with another method,
G\"{o}ttsche [G\"{o}] also considered the difference of Hodge numbers (and
hence Betti numbers) of moduli spaces of rank 2.
Matsuki and Wentworth [M-W] also generalized the chamber structure of
polarizations.
Combining another chamber structure,
they showed that the rational map between two moduli spaces is factorized to a
sequence of flips.
In sections 4 and 5, we assume that $X$ is a ruled surface which is not
rational.
Then, in the same way as in [Q1], we can give a condition for the existence of
stable sheaves.
Since we had computed the Picard group $Pic(\overline{M}_H(r;c_1,c_2))$ in case
of $(c_1,f)=0$ [Y3],
we assume that $0<(c_1,f)<r$.
In section 5, we shall compute the Picard group of $\overline{M}_H(r;c_1,c_2)$,
which is a generalization of [Q1] to $r \geq 3$.
The proof is the same as that in [D-N].
As is well known, it is difficult to treat the moduli spaces on rational ruled
surfaces (cf. [D-L], [R]).
However we can also check that $M_H(r;c_1,c_2)$ is emply or not in principle.

I would like to thank Professor S. Mori for valuable suggestions.

\section{chamber structure}

\subsection{}Notation

Let $X$ be a smooth projective surface defined over $\Bbb C$.
Let $NS(X)$ be the Neron-Severi group of $X$ and $Num(X)=NS(X)/\text{torsion}$.
Let $C(X) \subset Num(X) \otimes _{\Bbb Z}\Bbb R$ be the ample cone.
We denote the moduli space of stable sheaves of rank $r$ with Chern classes
$(c_1,c_2) \in H^2(X,\Bbb Q)\times H^4(X,\Bbb Q)$ by $M_H(r;c_1,c_2)$
and the
Gieseker-Maruyama compactification of $M_H(r;c_1,c_2)$ by
$\overline{M}_H(r;c_1,c_2)$.
We denote the open subscheme of $M_H(r;c_1,c_2)$ consisting of $\mu$-stable
sheaves
by $M_H(r;c_1,c_2)^{\mu}$ and the open subscheme consisting of $\mu$-stable
vector bundles by $M_H(r;c_1,c_2)^{\mu}_0$.
For a torsion free sheaf $E$ on $X$,
we set $\mu(E)=\frac{c_1(E)}{\rk(E)} \in H^2(X,\Bbb Q)$
and $\Delta(E)=\frac{1}{\rk(E)}(c_2(E)-\frac{\rk(E)-1}{2\rk(E)}c_1(E)^2) \in
H^4(X,\Bbb Q)$.
For a $x \in H^2(X,\Bbb Q)$, we set $P(x)=(x,x-K_X)/2+\chi(\cal O_X)$.
For a scheme $S$, we denote the projection $X \times S \to S$ by $p_S$.

\subsection{}

In this section, we shall generalize the chamber structure of polarizations in
[Q2].
For a torsion free sheaf $E$, we set $\gamma(E):=(\rk(E),\mu(E),\Delta(E)) \in
H^0(X,\Bbb Q) \times H^2(X,\Bbb Q) \times H^4(X,\Bbb Q)$.
For $\gamma \in \prod_{i=0}^2 H^{2i}(X, \Bbb Q)$,
 let $M_H^{\gamma}$ be the set of torsion free sheaves $E$ defined over $\Bbb
C$
with $\gamma(E)=\gamma$
which is $\mu$-semi-stable with respect to $H$.

\begin{lem}\label{lem:1}
Let $E$ be a torsion free sheaf which is defined by an extension
$ 0 \to F_1 \to E \to F_2 \to 0$. Then
$\Delta(E)=\frac{\rk(F_1)}{\rk(E)}\Delta(F_1)+\frac{\rk(F_2)}{\rk(E)}\Delta(F_2)-
\frac{\rk(F_1)\rk(F_2)}{2\rk(E)^2}((\mu(F_1)-\mu(F_2))^2)$.
\end{lem}

\begin{lem}
Let $B$ be a subset of $C(X)$.
 Let $\cal F_B(\gamma)$ be the set of filtrations $F:0 \subset F_1 \subset F_2
\subset \cdots \subset F_{s-1} \subset F_s=E$
 which satisfies (1) $\gamma(E)=\gamma$, (2) $\Delta_i=\Delta(F_i/F_{i-1}) \geq
0$ and (3)
there is an element $H \in B$ with $(\mu(F_{i-1})-\mu(F_i),H)=0$ for $2 \leq i
\leq s$.
If $B$ is compact, then $S_B(\gamma)=\{(\gamma(F_1),\cdots,\gamma(F_s))|
\text{$F_i$ is the $i$-th filter of $F \in \cal F_B(\gamma)$}\}$ is a finite
set.
\end{lem}

\begin{pf}
We denote $\mathrm{gr}_i(F):=F_i/F_{i-1}$ by $E_i$.
By using Lemma \ref{lem:1} successively, we see that
\begin{equation}\label {eq:2.1}
 \Delta(E)=\sum_{i=1}^s
\frac{\rk(E_i)}{\rk(E)}\Delta(E_i)-\sum_{i=2}^s\frac{\rk(E_{i-1})}{2 \rk(E_i)
\rk(E)}((\mu(F_{i-1})-\mu(F_i))^2).
\end{equation}
By the Hodge index theorem, we get $-((\mu(F_{i-1})-\mu(F_i))^2) \geq 0$ and
$-((\mu(F_{i-1})-\mu(F_i))^2)=0$ if and only if $\mu(F_{i-1})-\mu(F_i)=0$.
 By [F-M, II, Lemma 1.4], the set of ${c_1(F_i)}$ is finite.
 Hence ${\Delta_i}$ is finite. Therefore $S_B(\gamma)$ is a finite set.
\end{pf}

\begin{rem}
For a filtration $F:0 \subset F_1 \subset F_2 \subset \cdots \subset F_{s-1}
\subset F_s=E$ which belongs to $\cal F_B(\gamma)$,
$F':0 \subset F_i \subset F_s$ belongs to $\cal F_B(\gamma)$ for $1 \leq i \leq
s-1$.
In fact \eqref{eq:2.1} implies that $\Delta(F_i)\geq 0$ and
$\Delta(F_s/F_i)\geq 0$.
\end{rem}

\begin{defn}
For an element $F:0 \subset F_1 \subset F_2 \subset \cdots \subset F_{s-1}
 \subset F_s=E$ of $\cal F_{C(X)}(\gamma)$, we define a wall $W^F:=\cup_{i}\{ H
\in C(X)|(\mu(F_s)-\mu(F_{i}),H)=0 \}$, where $i$ runs for $1 \leq i \leq s-1$
with $\mu(F_s)-\mu(F_i) \ne 0$.
 By the above lemma, $\cup_F W^F$ is locally finite.
We shall call the connected component of $C(X) \setminus \cup_F W^F$ by
chamber.
\end{defn}

\begin{lem}\label{lem:2}
 Let $H$ and $H'$ be ample divisors which belong to a chamber $\cal C$.
Let $E$ be a $\mu$-semi-stable sheaf with respect to $H$.
Then $E$ be also $\mu$-semi-stable with respect to $H'$,
and hence we may denote $M_H^{\gamma}$ by $M_{\cal C}^{\gamma}$.
\end{lem}

\begin{pf}
Assume that $E$ is not $\mu$-semi-stable with respect to $H'$.
We shall construct a wall which separates $H'$ from $H$.
There is a filtration $F$ of $E$ such that $(\mu(F_{i-1})-\mu(F_i),H')>0$,
$2\leq i \leq s$
and $\Delta(\gr_i(E))\geq 0$, $1 \leq i \leq s$.
In fact, let $F:0 \subset F_1 \subset F_2 \subset \cdots \subset F_{s-1}
 \subset F_s=E$ be the Harder-Narasimhan filtration of $E$ with respect to
$H'$.
Then, the Bogomolov-Gieseker inequality implies that $\Delta(F_i/F_{i-1}) \geq
0$.
Let $H_t=H'+t(H-H'), 0 \leq t \leq 1$ be a line segment joining $H$ and $H'$.
There is a $t_1 \in \Bbb Q$ such that $(\mu(E_i)-\mu(E_{i+1}),H_t)>0$
 for $t<t_1$ and $(\mu(E_j)-\mu(E_{j+1}),H_{t_1})=0$ for some $j$.
Let $\{F_1',F_2',\dots,F_{s(t_1)}' \}$ be a subset of $\{F_1,F_2,\dots,F_s \}$
such that
$(\rk(F_i'),(F_i',H_{t_1})) \in  \Bbb Q \times \Bbb Q$, $1 \leq i \leq s(t_1)$
are vertices of the convex hull of
$\{(\rk(F_i),(F_i,H_{t_1})) \}_{i=1}^{s(t_1)}$.
By using Lemma \ref{lem:1}, we see that $\Delta(F_i'/F_{i-1}') \geq 0$.
Assume that $s(t_1) \ne s$.
Applying this argument successively,
we obtain a filtration $F'':0 \subset F_1'' \subset F_2'' \subset \cdots
\subset F_u''=E$
such that $\Delta(F_i''/F_{i-1}'')\geq 0$ and
$(\mu(F_i'')-\mu(F_{i+1}''),H_{t'})=0$
for some $t'$ with $0 < t' \leq 1$, moreover $\mu(F_i'')-\mu(F_{i+1}'')\ne 0$.
This implies that $H_{t'}$ belongs to a wall, which is a contradiction.
\end{pf}

\begin{defn}
Let $W$ be a wall and $\cal C$ a chamber such that $\overline{\cal C}$
intersects $W$.
Let $H$ be an ample divisor belonging to $\overline{\cal C} \cap W$ and $H_1$
an ample divisor which belongs to $\cal C$.
$V_{H,\cal C}^{\gamma}$ be the set of $\mu$-semi-stable sheaves with respect to
$H$
such that $E$ is not $\mu$-semi-stable with respect to $H_1$ and
$\gamma(E)=\gamma$.
\end{defn}

We shall investigate the set $V_{H,\cal C}^{\gamma}$.
We set $H_t=H_1+t(H-H_1)$, and $B=\{H_t| 0 \leq t \leq 1 \}$.
For an element $E$ of $V_{\cal C,H}^{\gamma}$, $S_B(\gamma)$ is a finite set.
Hence $S=S_B(\gamma) \cup \union\limits_{F \in \cal F_B(\gamma)}\cup_i
S_B(\gamma(\gr_i(E)))$ is a finite set.
Then there is a number $t'$ such that for all $t$ with $t' \leq t < 1$,
$F$ is the Harder-Narasimhan filtration of $E$ with resrect to $H_t$ if and
only if $F$ is that with respect to $H_{t'}$.
In fact, let $W^G$ be a wall defined by a $(\gamma(G_1),\cdots,\gamma(G_s))\in
S$ and
$I=\{H_t|t' \leq t \leq 1\}$ an interval which is contained in $B \setminus
\cup_{G}W^G$.
Let $F:0 \subset F_1 \subset F_2 \subset \cdots \subset F_{s-1}
 \subset F_s=E$ be the Harder-Narasimhan filtration of $E$ with respect to
$H_{t'}$.
In the same way as in the proof of Lemma \ref{lem:2}, there is a subset
$\{F_1',F_2',\dots,F_{s'}' \}$ of $\{F_1,F_2,\dots,F_s \}$
such that $F':0 \subset F_1' \subset F_2' \subset \cdots \subset F_{s'}'=E$
belongs to $\cal F_I(\gamma)$.
By the choice of $t'$, we get $(\mu(F_i')-\mu(F_{i+1}'),H)=0$.
 Moreover since $S_I(\gamma(F'_i/F'_{i-1}))$ is a subset of $S$,
$\{F_1',F_2',\dots,F_{s'}' \}$ must be $\{F_1,F_2,\dots,F_s \}$.
Thus $F$ belongs to $\cal F_B(\gamma)$.
If $\gr_i(E)$ is not $\mu$-semi-stable with respect to some $H_t$ with $t'<t
\leq 1$,
then $t'$ and $t$ are separated by a wall (Lemma 2.3), which is a
contradiction.
Thus $F$ is the Harder-Narasimhan filtration of $E$ for all $H_t$, $t' \leq
t<1$.
Therefore we get the following proposition.

\begin{prop}
Let $C$ be a 2-dimensional vector space such that $C \cap \cal C \ne \emptyset$
and $H \in C$.

(1)There is an element $H_1 \in C$ and
$V_{H,\cal C}^{\gamma}$ is the set of torsion free sheaves $E$ such that
$E$ has the Harder-Narasimhan filtration $F$ with respect to $H_1$
which is also Harder-Narasimhan filtration with respect to $H_t, 0 \leq t<1$,
and $F$ belongs to $\cal F_{\{H\}}(\gamma)$.

(2) $M_H^{\gamma}=M_{\cal C}^{\gamma} \amalg V_{H,\cal C}^{\gamma}$.

\end{prop}

\section{Equivariant cohomology of $M_H^{\gamma}$ }

\subsection{}
Let $C$ be a smooth projective curve with genus $g$
and $\pi:X \to C$ a ruled surface.
Let $C_0$ be a minimal section of $\pi$ with $(C_0^2)=-e$.
We assume that $e>\chi=2g-2$.
In this section, we shall define a cohomology of $M_H^{\gamma}$ and consider
the effect of change of polarizations.
For a scheme $S$, we denote the projection $S \times X \to S$ by $p_S$.
 Let $D=nH$, $n \gg 0$ be an ample divisor such that for an element $E \in
M_H^{\gamma}$,
$E(D)$ is generated by global sections and $H^j(X,E(D))=0$ $j>0$.
 Let $Q^{\gamma}$ be an open subscheme of $Quot_{\cal O_X(-D)^{\oplus N}/X/\Bbb
C}$ such that
for a quotient $\cal O_X(-D)^{\oplus N} \to E$, $E$ belongs to $M_H^{\gamma}$,
$H^0(X, \cal O_X^{\oplus N}) \cong H^0(X, E(D))$ and $H^j(X, E(D))=0 ,j>0$.
Since $(K_X,H)<0$, $Q^{\gamma}$ is smooth.
Let $H^*_{GL(N)}(Q^{\gamma},\Bbb Q):=H^*(Q^{\gamma} \times_{GL(N)}E(GL(N)),\Bbb
Q)$
be the equivariant cohomology of $Q^{\gamma}$,
where $E(GL(N))$ is the universal $GL(N)$-bundle over the classifying space.

\begin{lem}
$H^*_{GL(N)}(Q^{\gamma},\Bbb Q)$ does not depend on the choice of $Q^{\gamma}$.
 We denote this cohomology by $\tilde H^*(M_H^{\gamma},\Bbb Q)$ and the
Poincar\'{e} polynomial $\sum_i \dim \tilde{H}^i(M_H^{\gamma},\Bbb Q) z^i$
by $\tilde P (M_H^{\gamma},z)$.

\end{lem}
\begin{pf}
 Let $Q_i^{\gamma}$ $(i=1,2)$ be an open subscheme of $Quot_{\cal
O_X(-D_i)^{\oplus N_i}/X/\Bbb C}$
which satisfies the above conditions.
 Let $q_i:\cal O_{Q_i^{\gamma} \times X}(-D_i)^{\oplus N_i} \to \cal U_i$
be the universal quotient on $Q_i^{\gamma} \times X$.
 From the construction, $p_{Q_1^{\gamma}*}\cal U_1(D_2)$ is a locally free
sheaf on $Q_1^{\gamma}$.
 Let $\varphi: \Bbb V=\Bbb  V(\cal Hom(\cal O_{Q_1^{\gamma}}^{\oplus N_2},
p_{Q_1^{\gamma}*}\cal U_1(D_2))^{\vee}) \to Q_1^{\gamma}$ be a vector bundle
over $Q_1^{\gamma}$
and $h_1:\cal O_{\Bbb V}^{\oplus N_2} \to \varphi^*p_{Q_1^{\gamma}*}\cal
U_1(D_2)$ the universal homomorphism.
 Let $\cal G$ be the open subscheme of $\Bbb V$ such that $h_1$ is an
isomorphism.
 Then $\cal G$ is a principal $GL(N_2)$-bundle over $Q_1^{\gamma}$ and
there is a surjection $\cal O_{\cal G \times X}(-D_2)^{\oplus N_2} \to \cal
U_1$.
 For a $S$-valued point of $\cal G$, there is a flat family of quotients
$\cal O_{S \times X}(-D_1)^{\oplus N_1} \to \cal E$
and an isomorphism $\cal O_S^{\oplus N_2} \cong p_{S*}\cal E(D_2)$.
 It defines a surjection $q:\cal O_{S \times X}(-D_2)^{\oplus N_2} \to \cal E$,
and conversely for a surjection $q$,
it defines an isomorphism $\cal O_S^{\oplus N_2} \cong p_{S*}\cal E(D_2)$.
 Thus we obtain the following.
\begin{equation}
\cal G(S)=\left.\left\{(\cal E,q_1,q_2)\left|
\begin{aligned}
& \text{$\cal E$ is a flat family of coherent sheaves which belong to
$M_H^{\gamma}$,}\\
& \text{ and $q_i:\cal O_{S \times X}(-D_i)^{\oplus N_i} \to \cal E$ is a
surjective homomorphism.}
\end{aligned}
\right. \right\}\right/\sim
\end{equation}
where $(\cal E,q_1,q_2) \sim (\cal E',q_1',q_2')$ if and only if
there is an isomorphism $\psi:\cal E \to \cal E'$ with $q_i'=\psi \circ q_i$.
Hence $\cal G$ can be regarded as a $GL(N_1)$-bundle over $Q_2^{\gamma}$.
For simplicity, we denote $GL(N_i)$ by $G_i$.
$\cal G$ has a natural $G_1\times G_2$-action and $\cal G \times_{G_1 \times
G_2} EG_1 \cong Q_1^{\gamma} \times_{G_1} EG_1$.
 Therefore $\cal G \times_{G_1 \times G_2}EG_1 \times EG_2 \to Q_1^{\gamma}
\times _{G_1}EG_1$ is a $EG_2$-bundle.
 Thus $H^*_{G_1\times G_2}(\cal G,\Bbb Q) \cong H^*_{G_1}(Q_1^{\gamma},\Bbb
Q)$.
In the same way $H^*_{G_1\times G_2}(\cal G,\Bbb Q) \cong
H^*_{G_2}(Q_2^{\gamma},\Bbb Q)$.
  Hence $H^*_{G_1}(Q_1^{\gamma},\Bbb Q)  \cong H^*_{G_2}(Q_2^{\gamma},\Bbb Q)$,
which implies that $\tilde{H}^*(M_H^{\gamma},\Bbb Q)$ is well defined.
\end{pf}

\subsection{}
Let $\cal C$ be a chamber and $W$ a wall with $\overline{\cal C} \cap W \ne
\emptyset $.
Let $H$ be an ample divisor on $W$ and $H' \in \cal C$ an ample divisor which
is sufficiently close to $H$.
For a sequence of $\gamma_i=(r_i,\mu_i,\Delta_i)$,
$1 \leq i \leq s$ with $(\mu_i-\mu_{i+1},H)=0$ and $(\mu_i-\mu_{i+1},H')>0$,
we set
\begin{equation}
V_{H,\cal C}^{\gamma_1,\cdots,\gamma_s}=\left\{E \left|
\begin{aligned}
& \text{$E$ is not $\mu$-semi-stable with respect to $H'$ and for the
Harder-}\\
& \text{Narasimhan filtration $F:0 \subset F_1 \subset \cdots \subset F_s=E$,
$\gamma(F_i/F_{i-1})= \gamma_i$.}
\end{aligned}
\right.
\right \}.
\end{equation}

Let $Q^{\gamma_1,\cdots,\gamma_s}$ be the subscheme of $Q^{\gamma}$ such that
for a quotient $\cal O_X(-D)^{\oplus N} \to E$, $E$ belongs to $V_{H,\cal
C}^{\gamma_1,\cdots,\gamma_s}$,
and $\Gamma_{H,\cal C}$ the set of sequence $(\gamma_1,\cdots,\gamma_s)$.
By [D-L, 1], $Q^{\gamma_1,\cdots,\gamma_s}$ is a smooth locally closed
subscheme of $Q^{\gamma}$.
For the same $D$, let $Q^{\gamma_i}$ be an open subscheme of $Quot_{\cal
O_X(-D)^{\oplus N_i}/X/ \Bbb C}$
a quotient $\cal O_X(-D)^{\oplus N_i} \to E$ is contained in $Q^{\gamma_i}$ if
and only if $E$ belongs to $M_H^{\gamma_i}$.
In the same way as in [Y2, Appendix] (cf. [A-B],[K]), we obtain the following
theorem.

\begin{thm}\label{thm:1}
(1) \begin{align*}
H_{GL(N)}^*(Q^{\gamma_1,\cdots,\gamma_s},\Bbb Q) & \cong
\otimes_iH_{GL(N_i)}^*(Q^{\gamma_i},\Bbb Q)\\
& \cong \otimes_i\tilde H^*(M_{\cal C}^{\gamma_i},\Bbb Q).
\end{align*}
\newline
(2)
$d_{\gamma_1,\cdots,\gamma_s}:=\codim Q^{\gamma_1,\cdots,\gamma_s}=
-\sum_{i<j}r_ir_j(P(\mu_j-\mu_i)-\Delta_i-\Delta_j)$.
\newline
(3)$$ \tilde P(M_H^{\gamma},z)=\tilde P(M_{\cal
C}^{\gamma},z)+\sum_{(\gamma_1,\cdots,\gamma_s)
\in \Gamma_{H,\cal C}}z^{2d_{\gamma_1,\cdots,\gamma_s}}\prod_{i=1}^s\tilde
P(M_{\cal C}^{\gamma_i}, z).$$
\end{thm}

\begin{pf}
The proof is the same as that in [Y2, Appendix], so we shall give a sketch of
the proof.
Let $q_i:\cal O_{Q^{\gamma_i} \times X}(-D)^{\oplus N_i} \to \cal F_i$ be the
universal quotient.
We set $Z=\prod_{i=1}^sQ^{\gamma_i}$ and denote the $i$-th projection by
$\varpi_i$.
Then the quotient $\oplus_i q_i:\oplus_{i=1}^s \varpi_i^* \cal O_{Q^{\gamma_i}
\times X}(-D)^{\oplus N_i} \to \oplus_{i=1}^s \varpi_i^* \cal F_i$
defines a morphism $Z \to Q^{\gamma}$, which is an immersion.
We set $Y_1=Z$.
We shall define a sequence of schemes $Y_s \to \cdots \to Y_2 \to Y_1$
and quotients $\oplus _{j=1}^i\cal O_{Y_j \times X}(-D)^{\oplus N_j} \to \cal
F_{1,2,\cdots,i}$
$1 \leq i \leq s$ as follows.
Let $\psi_2:Y_2 \to Y_1$ be the vector bundle defined by a locally free sheaf
$\Hom_{p_{Y_1}}(\varpi_2^*(\ker q_2), \varpi_1^* \cal F_1)$.
There is a family of quotients $q_{1,2}:\oplus_{i=1}^2 \cal O_{Y_2 \times
X}(-D)^{\oplus N_i} \to \cal F_{1,2}$, which induces $q_1$ and $q_2$.
For $\psi_i:Y_i \to Y_{i-1}$ and $\oplus _{j=1}^i\cal O_{Y_j \times
X}(-D)^{\oplus N_j} \to \cal F_{1,2,\cdots,i}$,
$\Hom_{p_{Y_i}}(\ker q_{i+1}, \cal F_{1,2,\cdots,i})$ is a locally free sheaf
on $Y_i$.
Let $q_{i+1}:Y_{i+1} \to Y_i$ be the associated vector bundle on $Y_i$.
Then there is a quotient $\oplus _{j=1}^{i+1}\cal O_{Y_j \times X}(-D)^{\oplus
N_j} \to \cal F_{1,2,\cdots,i+1}$.
Let $P_{\gamma_1, \cdots ,\gamma_s}$ be the parabolic subgroup of $GL(N)$ which
preserves the filtration
$0 \subset \cal O_{Z \times X}(-D)^{\oplus N_1} \subset \oplus _{i=1}^2  \cal
O_{Z \times X}(-D)^{\oplus N_i}
\subset \cdots \subset \oplus_{i=1}^s \cal O_{Z \times X}(-D)^{\oplus N_i}$.
Then $Q^{\gamma_1, \cdots ,\gamma_s} \cong GL(N) \times_{ P_{\gamma_1, \cdots
,\gamma_s}}Y_s$.
The assertions follow from this (cf. [A-B, 7]).
\end{pf}

\begin{cor}
Let $\cal C'$ be another chamber with $\overline{\cal C'} \cap W \ne
\emptyset$.
Then
$$ \tilde P(M_{\cal C'}^{\gamma},z)=\tilde P(M_{\cal
C}^{\gamma},z)+\sum_{(\gamma_1,\cdots,\gamma_s)
\in \Gamma_{H,\cal
C}}\left\{z^{2d_{\gamma_1,\cdots,\gamma_s}}\prod_{i=1}^s\tilde P(M_{\cal
C}^{\gamma_i}, z)
-z^{2d_{\gamma_s,\gamma_{s-1},\cdots,\gamma_1}}\prod_{i=1}^s\tilde P(M_{\cal
C'}^{\gamma_i}, z) \right \}.$$

\end{cor}

\begin{rem}
 In the same way, we denote the set of $\mu$-semi-stable sheaves defined over
$\Bbb F_q$
by $M_{\cal C}^{\gamma}(\Bbb F_q)$. By using [D-R], we see that
\begin{multline}\label{eq:7}
\sum_{E \in M_{\cal C'}^{\gamma}(\Bbb F_q)}\frac{1}{\# \Aut(E)}=
\sum_{E \in M_{\cal C}^{\gamma}(\Bbb F_q)}\frac{1}{\# \Aut(E)}\\
+\sum_{(\gamma_1,\cdots,\gamma_s)
\in \Gamma_{H,\cal C}}\left\{q^{d_{\gamma_s,\gamma_{s-1},\cdots,\gamma_1}}
\prod_{i=1}^s \sum_{E \in M_{\cal C}^{\gamma_i}(\Bbb F_q)}\frac{1}{\#\Aut(E)}
-q^{d_{\gamma_1,\cdots,\gamma_s}}\prod_{i=1}^s \sum_{E \in M_{\cal
C'}^{\gamma_i}(\Bbb F_q)}\frac{1}{\#\Aut(E)}
 \right \}.
\end{multline}
By using the Weil conjectures [De] and results of Kirwan [K], we can also
obtain this corollary
(cf. [Y2, Proposition 4.3]).
By using \eqref{eq:7}, [Y2, Theorem 0.1] and analoguous argument to the proof
of [Y1, Proposition 0.3], in principle, we can compute the Betti numbers of the
moduli spaces of stable sheaves of rank 3 on $\Bbb P^2$
(in case of $c_1=0$, see [Y2, 4]). For example, we obtain the following.
\begin{align*}
P(M_H(3;1,2),z) &=1+z^2+z^4,\\
P(M_H(3;1,3),z) &=1+2z^2+5z^4+8z^6+10z^8+8z^{10}+5z^{12}+2z^{14}+z^{16},    \\
P(M_H(3;1,4),z) &=1+2z^2+6z^4+12z^6+24z^8+38z^{10}+54z^{12}+59z^{14}\\
                &\qquad
\qquad+54z^{16}+38z^{18}+24z^{20}+12z^{22}+6z^{24}+2z^{26}+z^{28}.
\end{align*}

\end{rem}

\begin{rem}
Let $X$ be a K3 or an Abelian surface and assume that Bogomolov-Gieseker
inequality holds.
Then \eqref{eq:7} also holds.
By using induction on $r$, we see that $\sum\limits_{E \in M^{\gamma}_{\cal
C}(\Bbb F_q)}\frac{1}{\#\Aut(E)}$
does not depend on $\cal C$ (cf. [G\"{o}]).
\end{rem}

\section{The existence of stable sheaves.}

\subsection{}

In this section, we assume that $X$ is not rational, that is, $g \geq 1$ and
assume that $e>-\chi=2g-2$.
We denote $C_0+xf$ by $H_x$.
Let $W_x$ be a wall containing $H_x$ and let $\cal C_x^+$ (resp. $\cal C_x^-$)
a chamber
containing $H_x+\epsilon f$ (resp. $H_x-\epsilon{}f$) with $0< \epsilon \ll 1$.
For $(\gamma{}_1,\cdots,\gamma{}_s) \in \Gamma{}_{H_x,\cal C_x^-}$,
we shall prove that $d_{\gamma{}_1,\cdots,\gamma{}_s} \geq 2$.
Since $(\mu{}_j-\mu{}_i)^2<0$ and $\Delta{}_i \geq 0$,
it is enough to prove that $r_ir_j(\mu{}_j-\mu{}_i,K_X/2) \geq 1$ for $i <j$.
We denote $r_ir_j(\mu{}_j-\mu{}_i)$ by $aC_0-bf$, and then $a$ and $b$ are
positive integer.
A simple calculation shows that $r_ir_j(\mu{}_j-\mu{}_i,K_X/2)=b+(g-1+e/2)a
\geq 3/2$.
Therefore $d_{\gamma{}_1,\cdots,\gamma{}_s} \geq 2$.
In particular, if $M_{\cal C_x^+}^{\gamma{}}$ is not empty,
then $M_{\cal C_x^-}^{\gamma{}}$ is not empty.
{}From this we obtain the following proposition.

\begin{prop}\label{prop:1}
For a triplet $\gamma=(r,\mu,\Delta)$ with $0<(\mu,f)<1$,
there exists a $\mu$-semi-stable sheaf $E$ of $\gamma(E)=\gamma$ with respect
to $H_x$
if and only if $x \leq \frac{e}{2}+\frac{r^2}{r_1r_2} \Delta$.
\end{prop}

\begin{pf}
Assume that $M_{H_x}^{\gamma}$ is not emty.
Since $(K_X+f,H)<0$, the deformation theory implies that there is a
$\mu$-semi-stable sheaf $E$
and an exact sequence
\begin{equation}
0 \to F_1(C_0) \to E \to F_2 \to 0,
\end{equation}
 where $F_1$ and $F_2$ are torsion free sheaves with
$(\mu(F_1),f)=(\mu(F_2),f)=0$.
We denote $\rk(F_i), \mu(F_i)$ and $\Delta(F_i)$ by $r_i,\mu_i$ and $\Delta_i$
respectively.
$(\mu(F_1(C_0))-\mu(F_2),H)=(\mu_1-\mu_2+C_0,C_0+xf)=(\mu_1-\mu_2,C_0)-e+x \leq
0$.
Thus $x \leq -(\mu_1-\mu_2,C_0)+e$.
On the other hand, $\Delta(E)=\frac{r_1}{r} \Delta(F_1(C_0))+\frac{r_2}{r}
\Delta_2
-\frac{r_1r_2}{2r^2}((\mu(F_1(C_0))-\mu(F_2))^2) \geq
-\frac{r_1r_2}{2r^2}(-e+2(\mu_1-\mu_2,C_0))$.
Thus $\Delta(E) \geq \frac{r_1r_2}{2r^2}(x-e/2)$, and hence $x \leq
\frac{e}{2}+\frac{r^2}{r_1r_2} \Delta$.
We shall next prove that the above condition is sufficient.
Let $E$ be a vector bundle defined by the following exact sequence.
$$
 0 \to F_1(C_0) \to E \to F_2 \to 0,
$$
where $F_1$ (resp. $F_2$) is the pull-back of a semi-stable vector bundle of
rank $r_1$ (resp. $r_2$) on $C$
with degree $d_1=r_1d+\frac{r_1^2-r_1}{2}e-c_2$ (resp. $d-d_1$).
Then $E$ is $\mu$-semi-stable with respect to
$H'=C_0+(\frac{e}{2}+\frac{r^2}{r_1r_2} \Delta(E))f$.
For general $H_x$, the claim follows immediately.
\end{pf}

\begin{cor}\label{cor:1}
 If $H$ is not on a wall, then there is a $\mu$-stable sheaf.
Moreover, $\codim(\overline{M}_H(r;c_1,c_2) \setminus M_H(r;c_1,c_2)^{\mu})
\geq 2$.
\end{cor}

\begin{pf}
We may assume that $r \geq 4$.
Let $F:0 \subset F_1 \subset F_2 \subset \cdots \subset F_s=E$
be a Jordan-H\"{o}lder filtration of a $\mu$-semi-stable sheaf $E$.
We set $E_i=F_i/F_{i-1}$.
Then $-\chi(E_i,E_j)=\rk(E_i)\rk(E_j)(g-1+\Delta(E_i)+\Delta(E_j))$.
By Proposition \ref{prop:1}, we see that $\rk(E_i)\Delta(E_i) > e/4$.
Hence $-\chi(E_i,E_j)>(\rk(E_i)+\rk(E_j))e/4\geq e \geq 1$.
The claim follows from this, (cf. [D-L] and [Y3, Proposition 2.3]).
\end{pf}

\begin{cor}\label{cor:2}
 If $H$ is not on a wall, then $\overline{M}_H(r,c_1,c_2)$ is locally
factorial.
\end{cor}
\begin{pf}
 This follows from the above corollary and [D2].
\end{pf}
 \begin{rem}\label{rem:2}
For $(\gamma_1,\cdots ,\gamma_s) \in \Gamma_{H_x,\cal C_x^-}$,
$d_{\gamma_1,\cdots ,\gamma_s} \geq 3$.
To prove this assertion, we may assume that $s=2$.
In the same way, we denote $r_1r_2(\mu_2-\mu_1)$ by $aC_0-bf$.
Since $r_1r_2(\mu_2-\mu_1,K_X/2)=b+((g-1)+e/2)a$, the assertion holds unless
$a=b=1$.
Assume that $a=b=1$, and then $(\mu_1,f)$ or $(\mu_2,f)$ is not an integer.
Hence we may assume that $(\mu_1,f)$ is not an integer.
By Proposition \ref{prop:1}, $r_1 \Delta_1 \geq \frac{1}{2}(x-\frac{e}{2})$.
Since $x=b/a+e=e+1$, $r_1 \Delta_1 >1/2$.
Therefore, we get $d_{\gamma_1,\gamma_2} \geq 3$.
\end{rem}

\section{The Picard group of $\overline{M}_H(r;c_1,c_2)$}

\subsection{}
In this section, we shall compute the Picard group of
$\overline{M}_{H_x}(r;c_1,c_2)$
under the assumption that $H_x$ does not lie on a wall.
By Corollary \ref{cor:1} and \ref{cor:2},
$Pic(\overline{M}_{H_x}(r;c_1,c_2))=Pic(M_{H_x}(r;c_1,c_2)^{\mu})$.
By using [Y1, Theorem 0.4], we see that $\codim(M_{H_x}(r;c_1,c_2)^{\mu}
\setminus M_{H_x}(r;c_1,c_2)^{\mu}_0) \geq r-1$,
and hence we shall compute $Pic(M_{H_x}(r;c_1,c_2)^{\mu}_0)$.
For a $\mu$-stable vector bundle $E \in M_{H_x}(r;c_1,c_2)^{\mu}_0$, $E^{\vee}$
is $\mu$-stable.
Hence $Pic(M_{H_x}(r;c_1,c_2)^{\mu}_0) \cong Pic(M_{H_x}(r;-c_1,c_2)^{\mu}_0)$.
Therefore we may assume that $0<(r,f) \leq r/2$.

For $c_1=r_1C_0+df$, we set $r_2=r-r_1,d_1=r_1d+\frac{r_1^2-r_1}{2}e-c_2$ and
$d_2=d-d_1$.
We shall first define a morphism $\overline{M}(r;c_1,c_2) \to J^{d_1} \times
J^{d_2}$.
Let $\cal Q$ be an open subscheme of quot-scheme $Quot_{V/X/\Bbb C}$ such that
$\overline{M}_{H_x}(r;c_1,c_2) =\cal Q/PGL(N)$
and $V \otimes \cal O_{\cal Q \times X} \to \cal U$ the universal quotients.
$\cal L=\det((1_{\cal Q} \times \pi)_! \cal U(-C_0))$ is a line bundle on $\cal
Q \times C$.
It defines a morphism $\lambda_{\cal Q}:\cal Q \to J^{d_1}$.
It is easy to see that $\lambda_{\cal Q}$ is $PGL(N)$-invariant.
Thus we get a morphism $\lambda:\overline{M}(r;c_1,c_2) \to J^{d_1}$.
The line bundle $\det \cal U \otimes \cal O(-r_1C_0)\otimes \cal L^{\vee}$
defines a morphism
$\nu: \overline{M}(r;c_1,c_2) \to J^{d_2}$.
Therefore we obtain the required morphism $\lambda \times
\nu:\overline{M}(r;c_1,c_2) \to J^{d_1} \times J^{d_2}$.

For simplicity, we denote $M_H(r;c_1,c_2)_0^{\mu}$ by $M_0^{\mu}$.
We set
$$
M^0=\{E|E \in M_0^{\mu} \text{ and }E_{|\pi{}^{-1}(P)} \cong \cal
O_{\pi{}^{-1}(P)}(1)^{\oplus r_1}
\oplus \cal O_{\pi{}^{-1}(P)}^{\oplus r_2} \text{ for all $P \in C$} \}.
$$
Assume that $r_1 \ne 1$. Since $e >\chi$, we get $(K_X+f,H_x)<0$.
Then, by the deformation theory, we see that
$\codim(M_0^{\mu} \setminus M^0)\geq 2$.
Next we assume that $r_1=1$.
For a fibre $l$, we set
$$
Z_l=\{E|E \in M_0^{\mu},\;E_l \cong \cal O_l(1)^{\oplus 2}
\oplus \cal O_l(-1) \oplus \cal O_l^{\oplus r_2-2} \}.
$$
Then we see that $\codim Z_l=2$ unless $Z_l=\emptyset$.
We set $Z=\cup_{P \in C}Z_{\pi{}^{-1}(P)}$. $Z$ is a locally closed subscheme
of $M_0^{\mu}$.
For a point $E$ of $Z$, there is an exact sequence
$0 \to L(C_0) \to E \to F \to 0$
, where $L$ is a line bundle with $c_1(L)=(d_1+1)f$ and $F$ is a torsion free
sheaf with
$c_1(F)=(d_2-1)f$ and $c_2(F)=1$.
We set $x_0=\frac{e}{2}+\frac{r^2}{r_1r_2} \Delta$ and
$x_1=\frac{e}{2}+\frac{r^2}{r_1r_2} (\Delta-\frac{1}{r})$.
If $x_1<x<x_0$, then $(H_x,\mu{}(L(C_0)))>(H_x,\mu{}(E))$, which is a
contradiction.
Thus $Z=\emptyset$.
Assume that $x<x_1$. Then for some $L$ and $F$ such that $F$ is semi-stable,
there is an exact sequence
\begin{equation}\label{eq:2}
0 \to L(C_0) \to E \to F \to 0
\end{equation}
such that $E$ is semi-stable (see the proof of Proposition \ref{prop:1}).
Thus $Z$ is not empty.
Therefore, to compute the Picard group of $\overline{M}_{H_x}(r;c_1,c_2)$, it
is enough to consider
$Pic(M^0 \cup (Z^0 \cap M_0^{\mu}) )$, where $Z^0$ is the subscheme of
$M_H(r;c_1,c_2)$
consisting of stable sheaves which are defined by the exact sequence
\eqref{eq:2}.

\subsection{}

We set $V_i=\cal O_X(-nH_x)^{\oplus N_i}$, $(i=1,2$).
Let $Quot_{V_1/X/\Bbb C}^{\gamma_1}$ ( resp. $Quot_{V_2/X/\Bbb C}^{\gamma _2}$)
be a quot-scheme parametrizing all quotients $V_1\to F_1$ (resp. $V_2 \to F_2$)
such that $\gamma(F_1)=(r_1,C_0+\frac{d_1}{r_1}f,0)$ (resp.
$\gamma(F_2)=(r_2,\frac{d_2}{r_2}f,0)$).
Let $Q_i$ ($i=1,2$) be the  open subscheme of $Quot_{V_i/X/\Bbb C}^{\gamma_i}$
consisting all quotients $V_i \to F_i$ which satisfy
\begin{enumerate}
\item $F_i$ is $\mu$-semi-stable with respect to $H_x$,
\item $F_{i|\pi^{-1}(\eta)}$ is a semi-stable vector bundle, where $\eta$ is
the generic point of $C$,
\item $H^0(X,V_i(nH_x)) \cong H^0(X,F_i(nH_x))$, $H^j(X,F_i(nH_x))=0,j>0$.
\end{enumerate}
Let $V_i \otimes \cal O_{Q_i \times X} \to \cal F_i$ be the universal quotient,
and $\cal K_i$ the universal subsheaf.
If we choose a sufficiently large integer $n$ and a suitable $N_i$,
then all $\mu$-semi-stable sheaves which satisfy (ii) are parametrized by
$Q_i$.
We set $\cal F_1'=\cal F(-C_0)$ and $\cal F_2'=\cal F_2$.
Let $g_i:G_i= Gr(p_{Q_i*}(\cal F_i'),r_i-1) \to Q_i$ be the grassmannian bundle
over $Q_i$
 parametrizing rank $r_i-1$ subbundle of $p_{Q_i*}(\cal F_i')$
   and $\cal U_i$ the universal subbundle of rank $r_i-1$.
   Since $p_{Q_i*}(\cal F_i')$ is $PGL(N_i)$-linearized,
   $G_i$ is $PGL(N_i)$-linearized.
  Let $G_i'=\{x \in G_i|\cal U_x \otimes \cal O_X \text{ is a subbundle of
}(\cal F_i')_x \}$.
  Let $h_i:D_i=\Bbb P(\cal Hom (\cal O_{G'_i}^{(r_i-1)},\cal U)^{\vee})\to
G'_i$ be a projective bundle
and $\cal O_{D_i}(1)$ the tautological line bundle on $D_i.$
On $D_i$, there is a homomorphism $\delta:\cal O_{D_i}^{\oplus (r_i-1)} \to
h_i^*\cal U \otimes \cal O_{D_i}(1)$.
Let $D_i'=\{x \in D_i |\delta_x \text{ is an isomorphism} \}$ be an open set of
$D_i$.
 Setting $\widetilde{\cal F_i'}=(g_i \circ h_i \times 1_X)^*\cal F_i'$ and
$\widetilde{\cal U}=p^*_{D_i'} h_i^* \cal U$,
there is an injective homomorphism on $D_i' \times X$: $\cal O_{D_i' \times
X}^{\oplus (r_i-1)} \to \widetilde{\cal U} \otimes p_{D_i'}^*\cal O_{D_i'}(1)
\to \widetilde{\cal F_i'} \otimes p_{D_i'}^*\cal O_{D_i'}(1)$.
 The quotient $\widetilde{\cal F_i'}\otimes p_{D_i'}^*\cal O_{D_i'}(1)/ \cal
O^{\oplus (r_i-1)}_{D_i' \times X}$
is a flat family of line bundles of degree $d_i$.
 Thus we obtain an extension
 \begin{equation}
 0 \to \cal O^{\oplus (r_i-1)}_{D_i' \times X} \to \widetilde{\cal F_i'}\otimes
p_{D_i'}^*\cal O_{D_i'}(1)
\to \det(\widetilde{\cal F_i'}\otimes p^*_{D_i'}\cal O_{D_i'}(1)) \to 0.
 \label{eq:1}
 \end{equation}
We set $Q=Q_1 \times Q_2$, $D=D_1' \times D_2'$ and $I=PGL(N_1)\times
PGL(N_2)$.
Then, in the same way as in [D-N, 7.3.4],
 we obtain the following exact sequence:
  \begin{equation}
   0 \to Pic^{I}(Q_1 \times Q_2) \to Pic^{I}(D_1 \times D_2) \to T \to 0,
   \label{eq:5}
   \end{equation}
  where $T$ is a finite abelian group with
$\#T=\frac{(r_1-1)d_1}{n_1}\frac{(r_2-1)d_2}{n_2}$.

Let $\cal P_i$ be a poincar\'{e} line bundle of degree $d_i$ on $J^{d_i}\times
X$.
 Let $\cal V_i:=\Ext^1_{p_{J^{d_i}}}(\cal P_i,\cal O_{J^{d_i} \times X}^{\oplus
(r_i-1)})$ be the relative extension sheaf on $J^{d_i}$.
The base change theorem implies that $\cal V_i$ is locally free.
Let $\mu_i:\Bbb P_i=\Bbb P(\cal V_i^{\vee}) \to J^{d_i}$ be the projection and
$\cal O_{\Bbb P_i}(1)$ the tautological line bundle on $\Bbb P_i$.

On $\Bbb P_i$, there is a universal family of extensions:

\begin{equation}
 0 \to \cal O_{\Bbb P_i \times X}^{\oplus (r_i-1)} \to \cal E_i \to \mu_i^*
\cal P \otimes \cal O_{\Bbb P_i}(-1) \to 0.
 \label{eq:4}
 \end{equation}
 We set $\Bbb P_i^{ss}=\{y \in \Bbb P_i|\text{$(\cal E_i)_y$ is
semi-stable}\}$.
 The extension \eqref{eq:1} gives a morphism $k':D_i' \to \Bbb P_i^{ss}$ such
that the pull--back of \eqref{eq:4} is \eqref{eq:1}.
Since $k'$ is $PGL(N_i)$-invariant and each fibre of $k'$ is an orbit of
$PGL(N_i)$,
[M-F-K, Proposition 0.2] implies that $\Bbb P_i^{ss}$ is a geometric quotient
of $D_i'$ by $PGL(N_i)$.
By [Y3, 4.2], $k'$ has a local section. Since the action of $PGL(N_i)$ is
set-theoretically free,
Z.M.T. implies that $D_i'$ is a Zariski locally trivial fibre bundle.
In particular, $Pic^{I}(D)=Pic(\Bbb P_1^{ss}\times \Bbb P_2^{ss})$ ([SGA I,
8]).

The base change theorem implies that
$\cal V=\Ext^1_{p_Q}(\cal F_2,\cal F_1)$ is a locally free sheaf on $Q$.
Let $\Bbb P_Q=\Bbb P(\cal V^{\vee}) \to Q$ be the projective bundle associated
to $\cal V^{\vee}$
and $\cal O_{\Bbb P_Q}(1)$ the tautological line bundle.
On $\Bbb P(\cal V^{\vee}) $, there is a universal extension

\begin{equation}
 0 \to \cal F_1 \to \cal E \to \cal F_2 \otimes \cal O_{\Bbb P_Q}(-1) \to 0.
 \label{eq:6}
 \end{equation}

 $I$ acts on $\Bbb P(\cal V^{\vee})$ and $\cal E$ is a $GL(N_1) \times
PGL(N_2)$-linearized.
 Let $\Bbb P_Q^s$ be the open subscheme of $\Bbb P(\cal V^{\vee})$
parametrizing stable sheaves.
 Then there is a surjective morphism $\lambda:\Bbb P_Q^s \to M^0$.
 It is easy to see that $\lambda$ is $I$-invariant and each fibre is an orbit
of this action.
 Thus $M^0$ is a geometric quotient of $\Bbb P_Q^s$ by $I$.
Let $S \to M^0$ be a smooth and surjective morphism such that there is a
universal family $\cal E$.
Then there is an exact sequence $0 \to \cal G_1(C_0) \to \cal E \to \cal G_2
\to 0$,
where $\cal G_1=\pi^*\pi_* \cal E(-C_0)$ and $\cal G_2=\cal E/\cal G_1(C_0)$.
There is an open covering $\{U_i \}$ of $S$ such that
$p_{U_i*} \cal G_1$ and $p_{U_i*} \cal G_2$ are free $\cal O_{U_i}$-module.
Then it defines a morphism $U_i \to Q$ and hence we get $U_i \to \Bbb P_Q^s$.
Therefore we get a local section of $\Bbb P_Q^s \times _{M^0}S \to S$.
Since $S \to M^0$ is a smooth morphism, $\Bbb P_Q^s \times _{M^0}S$ is smooth.
 By using Z.M.T., we see that $\Bbb P_Q^s \times _{M^0}S \to Q$ is a Zariski
locally trivial $I$-bundle.
By the descent theory ([SGA I, 8]), $\Bbb P_Q^s \to M^0$ is a $I$-bundle.
We shall prove that $\codim(\Bbb P_Q\setminus \Bbb P^s_Q) \geq 2$.
In the same way as in the proof of Theorem\ref{thm:1},
we define $Y_2=\Bbb V(\Hom_{p_Q}(\cal K_2,\cal F_1)^{\vee})$.
Then, there is a quotient $\cal O_{Y_2\times X}(-D)^{\oplus N} \to \cal
F_{1,2}$.
It defines a closed immersion $Y_2 \hookrightarrow \cal Q$.
Let $Y_2^s$ be the open subscheme of $Y_2$
parametrizing all quotients which are stable with respect to $H_x$.
Corollary \ref{cor:1} implies that $\codim(Y_2\setminus Y_2^s)\geq 2$.
Note that each fibre of $Y_2 \to \Bbb V(\Ext_{p_Q}(\cal F_2,\cal F_1)^{\vee})$
is contained in an orbit of $P_{\gamma_1,\gamma_2}$,
where $P_{\gamma_1,\gamma_2}$ is the parabolic subgroup of $GL(N)$
defined in the proof of Theorem \ref{thm:1}.
{}From this, we obtain that $\codim(\Bbb P_Q \setminus \Bbb P^s_Q) \geq 2$.

On $D_i$, there is a $GL(N_i)$-linearized line bundle such that
the action of the center $\Bbb C^{\times}$ is multiplication by constants.
In the same way, we obtain the following exact and commutative diagram:

\begin{equation}
 \begin{CD}
  @.0@.0@.0@. \\
  @.@VVV @VVV @VVV @. \\
  0 @>>> Pic^{I}(Q) @>>> Pic^{I}(\Bbb P_Q^s) @>>> \frac{n_1 n_2}{n} \Bbb Z @>>>
0 \\
  @. @VVV @VVV @VVV @. \\
  0  @>>> Pic^{I}(D) @>>> Pic^{I}(\Bbb P_D^s) @>>> \Bbb Z @>>> 0 \\
  @. @VVV @VVV @VVV @. \\
  0 @>>> T @>>> T' @>>> \Bbb Z/  \frac{n_1 n_2}{n}\Bbb Z  @>>> 0 \\
  @.@VVV @VVV @VVV  @. \\
  @.0 @.0@.0@.
 \end{CD}
\end{equation}
where $T'$ is a finite abelian group with
$\#T'=\frac{(r_1-1)d_1(r_2-1)d_2}{n}$.

\subsection{}

Let $K(X)$ be the Grothendieck group of $X$. Let $K^0(X)$ be the subgroup of
$K(X)$
which is generated by $\cal O_X-\cal O_X(-D)$ and $\cal O_C-\cal O_C(-D)$,
$D,D' \in Pic^0(X)$.
Then $K^0(X) \cong Pic^0(X)\oplus Alb(X)$.
We shall represent the class in $K(X)$ of $\cal O_X,\cal O_X(-f),
\cal O_X(-C_0)$ and $\cal O_X(-C_0-f)$ by $e_1,e_2,e_3$ and $e_4$ respectively.
 Then $K(X) \cong K^0(X) \oplus L$, where $L$ is the free $\Bbb Z$-module of
rank 4 generated by $e_i$, $1 \leq i \leq 4$.
  Let $\varepsilon$ be the class in $K(X)$ of a torsion free sheaf of rank $r$
with Chern classes $c_1,c_2$
and let $K(r;c_1,c_2)$ be the kernel of a homomorphism $K(X) \to \Bbb Z:x
\mapsto \chi(\varepsilon \otimes x)$.
Let $\cal E$ be a family of stable sheaves of rank $r$ with Chern classes
$c_1,c_2$ parametrized by a smooth scheme $S$.
Then $\det(p_{S!}(\cal E \otimes x))$, $x \in K(r;c_1,c_2)$  defines a line
bundle on $S$.
Thus we obtain a homomorphism $\kappa{}_S:K(r;c_1,c_2) \to Pic(S)$.
We can also define $\kappa{}:K(r;c_1,c_2) \to Pic(M_H(r;c_1,c_2))$, (see [Y3,
4.3]).
$K(r;c_1,c_2)=K^0(X) \oplus K$ where $K=K(r;c_1,c_2) \cap L$.

\begin{lem}\label{lem:6}
If $r_1 \ne 1$, then $K(r;c_1,c_2) \to Pic(M_H(r;c_1,c_2))/Pic(J^{d_1} \times
J^{d_2})$ is surjective.
\end{lem}

\begin{pf}
We denote $Pic^I(\Bbb P_D)/Pic(J^{d_1} \times J^{d_2})$ by $N$.
Since $Pic^I(D)/Pic(J^{d_1} \times J^{d_2}) \cong
 Pic(\Bbb P_1 \times \Bbb P_2)/Pic(J^{d_1} \times J^{d_2}) \cong
 \Bbb Z^{\oplus 2}$,
we get $N \cong \Bbb Z^{\oplus 3}$.
We shall prove that $\#(N/\im(\kappa_{\Bbb P_D}))=\#T'$.
 We denote the image of $\cal O_{\Bbb P_1}(1), \cal O_{\Bbb P_2}(1)$ and $\cal
O_{\Bbb P_D}(1)$ to $N$ by $\nu_1,\nu_2$ and $\nu$ respectively.
Let $\theta:Pic^I(\Bbb P_D) \to N$ be the quotient homomorphism.
We shall define $A_i$ $(1 \leq i \leq 4)$ as follows:
\begin{equation}
\left\{
\begin{split}
A_1 &:=\theta (\det p_{\Bbb P_D!}\cal
E)=-((2\chi+2d_1-e)\nu_1+(\chi+d_2)\nu_2+(r_2\chi+d_2)\nu),\\
A_2 &:=\theta (\det p_{\Bbb P_D!}\cal
E(-f))=-((2\chi+2d_1-2-e)\nu_1+(\chi+d_2-1)\nu_2+(r_2\chi+d_2-r_2)\nu),\\
A_3 &:=\theta (\det p_{\Bbb P_D!}\cal E(-C_0))=-(\chi+d_1)\nu_1,\\
A_4 &:=\theta (\det p_{\Bbb P_D!}\cal E(-C_0-f))=-(\chi+d_1-1)\nu_1.
\end{split}
\right.
\end{equation}

Let $\phi:L \to N $ be a homomorphism such that $\phi(e_i)=A_i$.
Then a simple calculation shows that $\phi(K)=\im(\kappa_{\Bbb P_D})$.

There is the following exact and commutative diagram

\begin{equation}
 \begin{CD}
  @.@.0@.0@. \\
  @.@.@VVV @VVV @. \\
  @.@. \ker \phi @= \ker \phi @. \\
  @.@.@VVV @VVV @. \\
  0 @>>> K @>>> L @> \psi>> n\Bbb Z @>>> 0 \\
  @. @| @VV{\phi}V @VVV @. \\
  0 @>>> K @>>> \phi(L) @>>> \phi(L)/K @>>> 0 \\
  @.@.@VVV @VVV  @. \\
  @.@. 0 @.0@.
   \end{CD}
\end{equation}
It is easy to see that $\ker \phi$ is generated by
$(\chi+d_1-1)e_3-(\chi+d_1)e_4$ and $N/\phi(L) \cong \Bbb Z/(r_2-1)d_2 \Bbb Z$.
 Hence $\psi(\ker \phi)=(r_1-1)d_1 \Bbb Z$.
 Therefore $\#N/K=\#(N/\phi(L))\# (\phi(L)/K)=\frac{(r_1-1)d_1(r_2-1)d_2}{n}$,
and hence $\# N/K=\#T'$.
Thus, we obtain our lemma.
\end{pf}

  \begin{lem}\label{lem:14}
  The restriction of $\kappa:K(r;c_1,c_2) \to Pic(M(r;c_1,c_2))$ to $K^0(X)$ is
injective and its image is $(\lambda \times \det)^*(Pic^0(J^{d_1} \times
J^{d_2}))$.
  \end{lem}
  \begin{pf}
  If $D=\sum_ia_i \pi^*(R_i), a_i \in \Bbb Z$, then we see that $\kappa_{\Bbb
P_D}(\cal O_X(D)-\cal O_X)=
\otimes_i (\cal P_{1R_i}^{\otimes 2}\otimes \cal P_{2R_i})^{\otimes a_i}$ and
$\kappa_{\Bbb P_D}(\cal O_{C_0}(D)-\cal O_{C_0})
=\otimes_i \cal P_{1R_i}^{\otimes a_i}$.
 The assertion follows immediately from this.
  \end{pf}

We shall first consider the case of $g \geq 2$.
\begin{thm}
Assume that $g \geq 2$ and $H_x$ does not lie on a wall.

(1) If $r_1 \ne 1$, then $Pic(\overline{M}_{H_x}(r;c_1,c_2))
\cong Pic(J^{d_1} \times J^{d_2}) \oplus \Bbb Z^{\oplus 3}$.

(2) If $r_1=1$, then
\begin{equation*}
Pic(\overline{M}_{H_x}(r;c_1,c_2)) \cong
\left\{
\begin{aligned} &  Pic(J^{d_1} \times J^{d_2}) \oplus \Bbb Z^{\oplus 2},
\;x_1<x<x_0\\
 & Pic(J^{d_1} \times J^{d_2}) \oplus \Bbb Z^{\oplus 3}, \;x<x_1.
\end{aligned}
\right.
\end{equation*}

(3) $Pic(\overline{M}_{H_x}(r;c_1,c_2))/Pic(J^{d_1}\times J^{d_2})$ is
generated by the image of $\kappa$.
\end{thm}

\begin{pf}
In the same way as in the proof of [Y3, Lemma 4.2],
we see that $\codim(\Bbb P_i^{ss}) \geq 2$, ($i=1,2$) for a sufficiently large
$d$.
Then (1) follows immediately from Lemma \ref{lem:6} and Lemma \ref{lem:14}.
We shall treat the case $r_1=1$.
In the same way, we can define $Q_1,Q_2,D_1$ and $D_2$, where $D_1=Q_1$.
Moreover we can define a projective bundle $\Bbb P_D$, where $D=D_1 \times
D_2$.
Then it is easy to see that $Pic^I(\Bbb P_D)/Pic(J^{d_1} \times J^{d_2}) \cong
\Bbb Z^{\oplus 2}$
and $\# (Pic^I(\Bbb P_D)/Pic^I(\Bbb P_{Q_1 \times
Q_2}))=\frac{(r_2-1)d_2}{n_2}$.
We denote $Pic^I(\Bbb P_D)/Pic^I(\Bbb P_{Q_1 \times Q_2})$ by $N$.
In the same way as in Lemma \ref{lem:6}, we denote the image of $\cal O_{\Bbb
P_2}(1)$ and
$\cal O_{\Bbb P_D}(1)$ to $N$ by $\nu_2$ and $\nu$ respectively.
Let $\phi:L \to N$ be the homomorphism such that
$\phi(e_1)=-((\chi+d_2)\nu_2+(r_2 \chi+d_2)\nu),
\phi(e_2)=-((\chi+d_2-1)\nu_2+(r_2 \chi+d_2-r_2)\nu),\phi(e_3)=\phi(e_4)=0$.
Then $\phi(K)=\im(\kappa_{\Bbb P_D})$.
There is the following exact and commutative diagram:
\begin{equation}
 \begin{CD}
  @.0@.0@.0@. \\
  @.@VVV @VVV @VVV @. \\
  0 @>>> K \cap \ker \phi @>>> \ker \phi @>>> \psi(\ker \phi) @>>> 0\\
  @.@VVV @VVV @VVV @. \\
  0 @>>> K @>>> L @> \psi>> \Bbb Z @>>> 0 \\
  @. @VVV @VV{\phi}V @VVV @. \\
  0 @>>> \phi(K) @>>> \phi(L) @>>> \phi(L)/K @>>> 0 \\
  @.@VVV @VVV @VVV  @. \\
  @. 0@. 0 @.0@.
   \end{CD}
\end{equation}
Since $\ker \phi=\Bbb Z e_3 \oplus \Bbb Z e_4$, we get $\Bbb Z/\psi(\ker
\phi)=0$,
and hence $\phi(K)=\phi(L)$.
 A simple calculation shows that $\# N/\phi(L)=\frac{(r_2-1)d_2}{n_2}$,
therefore $K \to Pic(M^0)/Pic(J^{d_1} \times J^{d_2})$ is surjective.
{}From this we get $Pic(M_{\cal C_{x_0}^-}(r;c_1,c_2))=Pic(J^{d_1} \times
J^{d_2}) \oplus \Bbb Z^{\oplus 2}$,
and hence we obtain the assertion for $x_1<x<x_0$.
We shall next prove the claim for $x<x_1$.
It is sufficient to compute $Pic(M_{\cal C_{x_1}^-}(r;c_1,c_2))$.

We set $V_i=\cal O_X(-nH_{x_1})^{\oplus N_i}$, $(i=3,4$).
Let $Quot_{V_3/X/\Bbb C}^{\gamma_3}$ ( resp. $Quot_{V_4/X/\Bbb C}^{\gamma _4}$)
be a quot-scheme parametrizing all quotients $V_3\to F_3$ (resp. $V_4 \to F_4$)
such that $\gamma(F_3)=(r_1,C_0+\frac{d_1+1}{r_1}f,0)$ (resp.
$\gamma(F_4)=(r_2,\frac{d_2-1}{r_2}f,\frac{1}{r_2})$).
Let $Q_i$ ($i=3,4$) be the  open subscheme of $Quot_{V_i/X/\Bbb C}^{\gamma_i}$
consisting quotients $V_i \to F_i$ which satisfy
\begin{enumerate}
\item $F_i$ is $\mu$-semi-stable with respect to $H_x$,
\item $F_{i|\pi^{-1}(\eta)}$ is a semi-stable vector bundle,
\item $H^0(X,V_i(nH_{x_1})) \cong H^0(X,F_i(nH_{x_1}))$,
$H^j(X,F_i(nH_{x_1}))=0,j>0$.
\end{enumerate}
Let $V_i \otimes \cal O_{Q_i \times X} \to \cal F_i$ be the universal quotient,
and $\cal K_i$ the universal subsheaf.
If we choose a sufficiently large integer $n$ and a suitable $N_i$,
then all $\mu$-semi-stable sheaves which satisfy (ii) are parametrized by
$Q_i$.
We set $R=Q_3 \times Q_4$.
The base change theorem implies that $\cal W=\Ext^1_{p_R}(\cal F_4,\cal F_3)$
is a locally free sheaf on $R$.
Let $\Bbb P_R=\Bbb P(\cal W^{\vee}) \to R$ be the projective bundle associated
to $\cal W^{\vee}$
and $\cal O_{\Bbb P_R}(1)$ the tautological line bundle.
On $\Bbb P(\cal W^{\vee}) $, there is a universal extension

\begin{equation}
 0 \to \cal F_3 \to \cal E \to \cal F_4 \otimes \cal O_{\Bbb P_R}(-1) \to 0.
 \label{eq:3}
 \end{equation}
 $PGL(N_3)\times PGL(N_4)$ acts on $\Bbb P(\cal W^{\vee})$ and $\cal E$ is a
$GL(N_3) \times PGL(N_4)$-linearized.
 Let $\Bbb P_R^s$ be the open subscheme of $\Bbb P(\cal W^{\vee})$
 parametrizing stable sheaves.
 Then there is a surjective morphism $\lambda:\Bbb P_R^s \to Z^0 \subset
M(r;c_1,c_2)$.
 It is easy to see that $\lambda$ is $PGL(N_3)\times PGL(N_4)$-invariant and
each fibre is an orbit of this action.
In the notation of {\bf 5.1}, we denote the pull--back of $Z^0$ to $\cal Q$ by
$\cal Z^0$.
[D-L, 1] implies that $\cal Z^0$ is smooth.
Hence $Z^0$ is a geometric quotient of $\Bbb P_R^s$ by $PGL(N_3)\times
PGL(N_4)$.
On $\cal Z^0$, there is an exact sequence
$ 0 \to \cal G_1(C_0) \to \cal U_{|\cal Z^0 \times X} \to \cal G_2 \to 0$,
where $\cal G_1$ is a flat family of line bundles with $c_1=(d_1+1)f$ and
$\cal G_2$ a flat family of torsion free sheaves of rank $r-1$ with Chern
classes $((d_2-1)f,1)$.
Then the normal bundle $\cal O_{\cal Z^0}(\cal Z^0)$ is isomorphic to
$\Ext^1_{p_{\cal Z^0}}(\cal G_1(C_0),\cal G_2)$.
It is easy to see that the pull--back of $\Ext^1_{p_{\cal Z^0}}(\cal
G_1(C_0),\cal G_2)$
to $\cal Z^0 \times _{Z^0} \Bbb P_R^s$ is isomorphic to the pull--back of
$\Ext^1_{p_{\Bbb P^s_R}}(\cal F_3,\cal F_4 \otimes \cal O_{\Bbb P_R}(-1))$ to
$\cal Z^0 \times _{Z^0} \Bbb P_R^s$.
Since $Pic(\Bbb P_R^s) \to Pic(\cal Z^0 \times _{Z^0} \Bbb P_R^s)$ is
injective, the pull--back of $\cal O_{Z^0}(Z^0)$ to $\Bbb P_R^s$ is isomorphic
to $\Ext^1_{p_{\Bbb P^s_R}}(\cal F_3,\cal F_4 \otimes \cal O_{\Bbb P_R}(-1))$.
By virtue of Remark \ref{rem:2}, we get $\codim(\Bbb P_R \setminus \Bbb P_R^s)
\geq 2$.
Let $\lambda' \times  \det:\overline{M}(r_2;(d_2-1)f,1) \to C \times J^{d_2-1}$
be the morphism defined in [Y3].
Let $J^{d_1+1}\times C \times J^{d_2-1} \to J^{d_1} \times J^{d_2}$ be the
morphism sending $(L,P,L')$ to
$(L \otimes \cal O_C(-P),L' \otimes \cal O_C(P))$.
Then the composition $Z^0 \to \overline{M}(r_1;(d_1+1)f,0) \times
\overline{M}(r_2;(d_2-1)f,1) \to
J^{d_1+1}\times C \times J^{d_2-1} \to J^{d_1} \times J^{d_2}$ is the same as
the restriction of
$\lambda \times \nu$ to $Z^0$.
It is easy to see that $K \cap \ker \phi=\Bbb
Z((\chi+d_1-1)e_3-(\chi+d_1)e_4)$,
and hence $L:=\kappa((\chi+d_1-1)e_3-(\chi+d_1)e_4))$ can be written as $\cal
O(nZ^0)\otimes  \cal L$,
where $\cal L$ is the pull--back of a line bundle on $J^{d_1} \times J^{d_2}$.
We shall prove that $n=-1$.
Since $\det(\cal E(-C_0))=\det(\cal F_3(C_0))\otimes \det(\cal F_4 \otimes \cal
O_{\Bbb P_R}(-1))$,
the restriction of $\det(\cal E(-C_0))$ to a fibre of $\Bbb P_R \to R$ is
$\cal O(1)$. From this, we see that $n=-1$.
(3) follows from the proof of (1) and (2).
\end{pf}

\subsection{}
We shall treat the case of $g=1$.
We assume that $(r_1,d_1) \ne 1$ and $(r_2,d_2)=1$, and then there is an
integers $r_2'$ and $d_2'$ such that
$r_2d_2'-r_2'd_2=1$ and $0<r_2'<r_2$. We set $r_2''=r_2-r_2'$ and
$d_2''=d_2-d_2'$.
Let $W$ be the subset of $M^0$ whose element $E$ has the following filtration
$F:0 \subset F_1
\subset F_2 \subset F_3=E$ such that
\begin{enumerate}
 \item $\rk(F_1)=r_1$, $c_1(F_1(-C_0))=d_1f$ and $c_2(F_1(-C_0))=0$,
 \item $\rk(F_2/F_1)=r_2'$, $c_1(F_2/F_1)=d_2'f$ and $c_2(F_2/F_1)=0$,
 \item $\rk(F_3/F_2)=r_2''$, $c_1(F_3/F_2)=d_2''f$ and $c_2(F_3/F_2)=0$.
\end{enumerate}

If $\frac{r}{r_1}\left(\frac{d_2''}{r_2''}-\frac{d}{r} \right) < x <
\frac{d_2}{r_2}-\frac{d}{r}$,
then for a general element $E$ of $W$, the Harder-Narasimhan filtration is
$F':0 \subset F_2 \subset F_3=E$.
Let $Q_2'$ be an open subscheme of $Quot_{V_2/X/\Bbb C}$ whose point $y$
satisfies that $\cal F_y$
is semi-stable or the Harder-Narasimhan filtration of $\cal F_y$ is $0 \subset
G_1 \subset \cal F_y$,
where $G_1$ is a semi-stable vector bundle of rank $r_2'$ and degree $d_2'$.
We shall replace $Q_2$ by $Q_2'$, and construct $\Bbb P_D$ and $\Bbb P_D^s$.
Let $\widetilde{W}$ be the open subscheme of $\Bbb P_D \setminus \Bbb P_D^s$
whose point defines an element of $W$.
Then, there is an exact sequence
$ \Bbb Z \widetilde{W} \to Pic^I(\Bbb P_D) \to Pic^I(\Bbb P_D^s) \to 0.$
In the same way, we see that $\cal O_{\widetilde{W}}(\widetilde{W})$ is a
primitive element
of $Pic(\widetilde{W})$.
Note that $Pic^I(\Bbb P_D)$ is isomorphic to
$Pic(\overline{M}_{H_{x'}}(r,c_1,c_2))$,
$x'<\frac{r}{r_1}(\frac{d_2''}{r_2''}-\frac{d}{r})$.
Therefore $Pic^I(\Bbb P_D^s) \cong Pic(J^{d_1} \times J^{d_2}) \oplus \Bbb
Z^{\oplus 3}/\Bbb Z \widetilde{W}
\cong Pic(J^{d_1} \times J^{d_2})\oplus \Bbb Z^{\oplus 2}$.
In the same way, we obtain the following theorem

\begin{thm}
If $g=1$, then
$Pic(\overline{M}_H(r;c_1,c_2)) \cong Pic(J^{d_1} \times J^{d_2}) \oplus \Bbb
Z^{\oplus a}$, $a=1,2$ or $3$.
\end{thm}

% Start-of-Trailer

% Local Variables:
% TeX-trailer-start: "% Start-of-Trailer"
% TeX-header-end: "% End-of-Header"
% TeX-master: t
% End:


\begin{thebibliography}{[M-F-K]}





\bibitem[A-B]{A-B:1}
 Atiyah, M. F., Bott, R.
{\it The Yang-Mills equations over Riemann surfaces,}
 Philos. Trans. Roy. Soc. London Ser. A {\bf308} (1982), pp. 523--615
\bibitem[De]{De:1}
Deligne, P.,
{\it La conjecture de Weil. I,}
Inst. Hautes Etudes Sci. Publ.Math. {\bf  43} (1974) pp. 273--307
\bibitem[D1]{D:1}
Drezet, J.-M.,
{\it Groupe de Picard des vari\'{e}t\'{e}s de modules de faisceaux semi-stables
sur $\Bbb P^2(\Bbb C)$,}
Ann. Inst. Fourier {\bf 38} (1988), pp. 105--168
\bibitem[D2]{D:2}
Drezet, J.-M.,
{\it Points non factoriels des vari\'{e}t\'{e}s de modules de faisceaux
semi-stable sur une surface rationelle,}
J. reine angew. Math. {\bf 413} (1991), pp. 99--126
\bibitem[D-L]{D-L:1}
Drezet, J.-M., Le-Potier, J.,
{\it Fibr\'{e}s stables et fibr\'{e}s exceptionnels sur $\Bbb P^2$,}
Ann.\ scient.\ \'{E}c. Norm. Sup., $4^e$ s\'{e}rie, t. {\bf 18} (1985), pp.
193--244
\bibitem[D-N]{D-N:1}
Drezet, J.-M., Narasimhan, M. S.,
 {\it Groupe de Picard des vari\'{e}t\'{e}s de modules de fibr\'{e}s
semi-stables sur les courbes alg\'{e}griques,}
 Invent.\ math. {\bf 97} (1989), pp.\ 53--94
\bibitem[D-R]{D-R:1}
Desale, U. V., Ramanan, S.
{\it Poincar\'e polynomials of the variety of stable bundles,}
Math. Ann. {\bf 216} (1975) pp. 233--244
\bibitem[F-M]{F-M:1}
Friedman, R., Morgan, J.,
{\it On the diffeomorphism types of certain algebraic surface I,}
J. Differential Geometry {\bf 27} (1988), pp. 371--398
\bibitem[Gi]{Gi:1}
Gieseker, D.,
{\it On a theorem of Bogomolov on Chern classes of stable bundles,}
Amer. J. Math. {\bf 101} (1979), pp. 77--85
\bibitem[G\"{o}]{G:1}
G\"{o}ttsche, L.,
{\it Change of polarization and Hodge numbers of moduli spaces of torsion free
sheaves on surfaces,}
MPI 94-9 (1994)
\bibitem[SGA I]{SGA:1}
Grothendieck, G.,
{\it S\'{e}minaire de g\'{e}om\'{e}trie alg\'{e}brique I,}
Springer Lecture Notes {\bf 224} (1971)
\bibitem[K]{K:1}
 Kirwan, F.,
{\it Cohomology of quotients in symplectic and algebraic geometry,}
 Princeton Math. Notes 31 Princeton N.J. (1984)
\bibitem[Mr]{M:2}
 Maruyama, M.,
 {\it Moduli of stable sheaves II,}
 J. Math. Kyoto Univ. {\bf 18} (1978), pp. 557--614
\bibitem[M-F-K]{M-F-K:1}
Mumford, D., Fogarty, J., Kirwan, F.,
{\it Geometric Invariant Theory,}
Springer-Verlag
\bibitem[M-W]{M-W:1}
Matsuki, K., Wentworth, R.,
{\it Mumford-Thaddeus principle on the moduli space of vector bundles
on an algebraic surface,}
preprint
\bibitem[Q1]{Q:1}
Qin, Z. B.,
{\it Moduli of stable sheaves on ruled surfaces and their Picard groups,}
J. reine angew. Math. {\bf 433} (1992), pp. 201--219
\bibitem[Q2]{Q:2}
Qin, Z. B.,
{\it Equivalence classes of polarizations and moduli spaces of sheaves,}
J. Differential Geometry {\bf 37} (1993) pp. 397--413
\bibitem[R]{R:1}
Rudakov, A. N.,
{\it A description of Chern classes of semistable sheaves on a quadric
surface,}
J. reine angew. Math. to appear
\bibitem[S]{S:1}
Str\o mme, S. A.,
{\it Ample Divisors on Fine Moduli spaces on the Projective Plane,}
 Math. Z. {\bf 187} (1984), pp. 405--423
\bibitem[Y1]{Y:1}
Yoshioka, K.,
{\it The Betti numbers of the moduli space of stable sheaves of rank 2 on $\Bbb
P^2$,}
J. reine angew. Math. to appear
\bibitem[Y2]{Y:2}
Yoshioka, K.,
{\it The Betti numbers of the moduli space of stable sheaves of rank 2 on a
ruled surface,}
Math. Ann. to appear
\bibitem[Y3]{Y:3}
Yoshioka, K.,
{\it The Picard group of the moduli space of stable sheaves on a ruled
surface,}
preprint





\end{thebibliography}
\end{document}